\documentstyle[12pt,twoside]{amsart}
\title{On $-K^2$ for normal surface singularities}
\address{Department of Mathematics, Zhongshan University,
Guangzhou, Guangdong, China\\
\ \ Department of Mathematics, Tokyo Institute of Technology,
Oh-Okayama, Meguro, Tokyo 152, Japan}
\author{Hao Chen\ \  \& \ \ Shihoko Ishii}

\newcommand{\bZ}{{\Bbb Z}}

\newcommand{\bN}{{\Bbb N}}
\newcommand{\bC}{{\Bbb C}}
\newcommand{\bQ}{{\Bbb Q}}
\newcommand{\ba}{{\bold a}}
\newcommand{\bm}{{\bold m}}
\newcommand{\bc}{{\bold c}}
\newcommand{\bn}{{\bold n}}
\newcommand{\Supp}{{\text{Supp}}}

\newtheorem{thm}{Theorem}[section]

\newtheorem{lem}[thm]{Lemma}
\newtheorem{cor}[thm]{Corollary}
\newtheorem{prop}[thm]{Proposition}

\theoremstyle{definition}
\newtheorem{defn}[thm]{Definition}

\newtheorem{say}[thm]{}
\newtheorem{exmp}[thm]{Example}

\theoremstyle{remark}

\begin{document}
\bibliographystyle{amsplain}
\maketitle
%%%%%%%%%%%
\pagestyle{myheadings}
\markboth{\hfill  CHEN \& ISHII\hfill}{\hfill $-K^2$ FOR SURFACE
SINGULARITIES\hfill}
%%%%%%%%%%%
\begin{abstract}
  In this paper we show the lower bound of the set of non-zero 
  $-K^2$ 
  for normal surface singularities establishing  
  that this set has no accumulation
  points from above.
  We also prove that every accumulation point from below is
  a rational number
  and every positive integer is an accumulation point.
  Every rational number can be an accumulation point modulo $\bZ$.
  We determine all accumulation points 
  in $[0, 1]$.
  If we fix the value  $-K^2$, then the  values of $p_g$, $p_a$,
  mult, embdim  and the numerical indices are bounded,
  while the numbers of the exceptional curves are not bounded.
\end{abstract}

\setcounter{section}{-1}

\section{\bf Introduction}
  For a normal surface  singularity $(X, x)$ over $\bC$,
  we have two kinds of plurigenera $\gamma_m(X,x)$ and
  $\delta_m(X,x)$  which are defined by Kn\"oller \cite{Kn} and
  Watanabe \cite{Wt} respectively.
  Both plurigenera grow in order at most 2 and the coefficients
  of the term of degree 2 are rational numbers.
  The leading coefficient of $\gamma_m(X, x)$ is $-K^2/2$,
  where  $K$  is the numerical canonical divisor on the 
  minimal resolution (cf. \cite{Mo}),
  and that of $\delta_m(X,x)$  is $-P^2/2$ (cf. \cite{Wah}).
  It is well known that $-K^2=0$  if and only if the singularity 
  is a rational double point.
  In this paper,
  we study the set of the values of $-K^2$.
  The set $\{-P^2\}$  is studied by Ganter in \cite{Ga}.
  Her results are:
  if one fixes the numerical index $m$  of  singularities,
  then non-zero $-P^2$ has the lower bound $1/42m$
  and the set of  $-P^2$ has no accumulation points 
  from above,
  which is equivalent to the descending chain condition 
  (D.C.C. for short)
  because of the lower bound.
  Here one should note that there are accumulation points of $-P^2$
  from above,
  if one does not fix the numerical index (see \ref{quotient}).
  Our results on $-K^2$ are simpler.
  The discussions go well without fixing the numerical index.
  We prove that non-zero $-K^2$  has the lower bound  $1/3$
  and the set of $-K^2$  has D.C.C..
  Then we show that all accumulation points from below 
  are rational numbers and every positive integer is an accumulation
  point.
  There are  many accumulation points so that every rational
  number can be an accumulation point modulo $\bZ$.
  The accumulation points in $[0, 1]$ turn out to be
  $\{ 1, m/(m+1)|m\in \bN\}$.
  Lastly we fix the value of $-K^2$ and observe the behavior
  of some invariants
  of singularities.
  We can see easily that the numbers of the exceptional curves
  are not bounded,
  but  we prove that invariants $p_g$, $p_a$, mult, embdim,
  and the numerical indices  are bounded.

  For compact surfaces, 
  this kind of problem is posed by Koll\'ar in \cite{Ko}
  and answered by Alexeev in \cite{Al} and reproved by Alexeev-
Mori in \cite{AM}.
  Our result for singularities on the property D.C.C. is 
  similar to the compact case.

  We would like to express our hearty thanks to Professors
  Jonathan Wahl and Masataka Tomari for the helpful 
  discussions with the second author.
  We are also grateful to the referee for the ideas 
  which make the proofs of 1.2, 2.3 and 2.8 simple.

%%%%%%%%%%%%%%%%%% section 1 %%%%%%%%%%%%%%%%%%%%%%%

\section{\bf Lower bound and D.C.C. for $-K^2$}

In this paper all surfaces are defined over $\bC$.

\begin{defn}
  Let $A=\sum A_i$  be a compact reduced  divisor whose intersection 
  matrix $(A_i\cdot A_j)_{i,j}$ is negative definite on 
  a normal $\bQ$-factorial surface $Y$.
  One defines the numerical canonical divisor $K(A)$ 
  as a $\bQ$-divisor with the support on $A$ as 
  $$K(A)\cdot A_i=K_Y\cdot A_i$$
  for every $i$.
  In particular,
  if  $Y$  is the minimal resolution of a normal singularity
  $(X,x)$ and $A$  is the exceptional divisor,
  then $K(A)$  is called the numerical canonical divisor of the 
  singularity  $(X,x)$.
\end{defn}
  Note that $K(A)$  is well defined because of the negative 
  definiteness of the intersection matrix.

  First one proves a basic property on $-K(A)^2$.
\begin{prop}
\label{subgraph}
%**(subgraph)
  Let  $A'<A$  be compact connected reduced divisors 
  whose intersection
  matrices are negative definite on non-singular surface $Y$.
  Assume that there is no (-1)-curve on  $A$.
  Then
  $$-K(A')^2\leq -K(A)^2.$$
\end{prop}
\begin{pf}
  Note that $-K(A)$ and $-K(A')$ are effective $\bQ$-divisors
  supported on $A$ and $A'$  respectively.
  Hence  $(K(A)-K(A'))\cdot C\geq K(A)\cdot C \geq 0$ 
  for a curve  $C\subset \Supp A$  with  $C\not\subset \Supp A'$,
  because $A$ has no (-1)-curves.
  On the other hand, 
  by the definition, $(K(A)-K(A'))\cdot C=0$ 
  for a curve  $C \subset \Supp A'$.
  Therefore $K(A)-K(A')$  is nef on $\Supp A$.
  Hence $0 \geq (K(A)-K(A'))(K(A)+K(A'))=K(A)^2-K(A')^2$
\end{pf}  
\vskip.5truecm
  Now one proves a lemma which is used for the lower bound and
  also D.C.C. for $-K^2$.
\begin{lem}
\label{basic lemma}
%**(basic lemma)
  Let  $A=\sum _{i=1}^rA_i$  be a compact reduced  
  divisor on a non-singular 
  surface.
  Assume that the intersection matrix is negative definite and
  $A$  has no (-1)-curves.
  Then
$$-K(A)^2\geq \sum_{i=1}^r(K(A)\cdot A_i)^2/(-A_i^2).$$
\end{lem}
\begin{pf}
  Since it is sufficient to prove the inequality for each
  connected component,
  one may assume that  $A$  is connected.
  If  $A$ consists of (-2)-curves,
  then  the equality holds,
  because the both hand sides are zero.
  So one assumes that $A$  has a non-(-2)-curve.
  Then, denoting $K=K(A)=\sum m_iA_i$,
  it follows that   $m_i < 0$ for every $i$,
  because $K$  has non-negative intersections with  all  $A_i$.

  Therefore 
$$(\ref{basic lemma}.1)\ \ \ \ \ \ K\cdot A_i=(\sum m_jA_j)\cdot
  A_i
  =m_iA_i^2+\sum_{j\neq i} m_jA_j\cdot A_i\leq m_iA_i^2.$$
  Substituting  $K\cdot A_i/(-A_i^2)\leq -m_i$ into
  $-K^2=-K\cdot (\sum m_i A_i)= \sum (-m_i)K\cdot A_i$,
  one obtains the inequality.
\end{pf}
\vskip.5truecm
\begin{thm}
\label{lower bound}
%**(lower bound)
   Let  $K$  be the numerical canonical divisor of a normal
  surface singularity  $(X,x)$.
  If $(X,x)$  is not rational double,
  then 
$$-K^2\geq \frac{1}{3}$$
  and the equality holds if and only if $(X,x)$ is
  the cyclic quotient singularity $X_{3,1}$.
\end{thm}
\begin{pf}
   Let  $A=\sum _{i=1}^rA_i$  be the exceptional divisor of 
  the minimal resolution,
  then it satisfies the condition of Lemma\ref{basic lemma}.
  Let $A_i$  be non-(-2)-curve and $p$  its arithmetic genus.
  Then $$(K\cdot A_i)^2/(-A_i^2)=(2p-2)^2/(-A_i^2)+2(2p-2)-A_i^2
  \geq 1/3,$$
  and the equality holds if and only if $A_i^2=-3$ and $p=0$.
  In fact, if $p>0$, then the left hand side is $\geq 1$;
  on the other hand,
  if $p=0$, then the left hand side is $\geq 1$ 
  for $-A_i^2\geq 4$
  and is $=1/3$ for $-A_i^2=3$.
  By Lemma\ref{basic lemma}, 
  it follows  $-K(A)^2\geq 1/3$.
  Here assume the equality,
  then there is only one non-(-2)-curve $A_i$ which is (-3)-curve
  and 
  the equality should hold in (\ref{basic lemma}.1);
  $\sum _{j\neq i}m_j A_j\cdot A_i=0$.
  Since  $m_j<0$ for every $j$ and $A$  is connected,
  the last equlaity implies that there is no other component.
\end{pf}

  In order to prove D.C.C. for $-K^2$ one should prepare
  the following notion.

\begin{defn}
  Let $A=\sum _{i=1}^rE_i$  and 
  $A'=\sum_{i=1}^r E'_i+\sum_{j=1}^nF_j$  be  
  negative definite divisors with 
  no (-1)-curves on non-singular surfaces $Y$ and $Y'$ respectively.
  One says that $A'$ is a (-2)-insertion of $A$,
  or $A \to A'$  is a (-2)-insertion,
  if the following holds:

  $E_1, E_2, E'_1, E'_2$ are (-2)-curves
  and $\sum_{j=1}^nF_j$ is a chain
  of (-2)-curves
$$p_a(E_i)=p_a(E'_i) \ \ \ (\forall i=1,\ldots ,r)$$
$$E_1\cdot E_2=1,\ \ E'_1\cdot E'_2=0,\ \ E_i\cdot E_k=
E'_i\cdot E'_k, \ \ for (i,k)\neq (1,2)$$
$$F_1\cdot E'_1=F_n\cdot E'_2=1$$
$$F_j\cdot E'_i=0 \ \ (\forall j\neq 1, n,\ \  \forall i).$$
  In this case one calls $\sum_{j=1}^nF_j$ the insertion string
  and $E'_1$, $E'_2$ the props of the insertion string.
\end{defn}
\vskip.5truecm
\begin{lem}
\label{property}
%**(property)
  Let $A'=\sum_{i=1}^r E'_i+\sum_{j=1}^nF_j$  be 
  a (-2)-insertion of 
  $A=\sum _{i=1}^rE_i$.
  Let $E'_1$  and  $E'_2$  be the props of the insertion string.
  Let  $\varphi$  be the contraction morphism of $F_1,\ldots , F_n$
  and denote $\varphi_*A'$ and $\varphi_*E'_i$  by $A''$ and 
  $E_i'' $, respectively. 
  Then the following hold:

(i)  if the intersection matrix $M$ of $A$ is
$$M=\pmatrix
   -2&1&^t\ba_1\\
   1&-2&^t\ba_2\\
     \ba_1& \ba_2&N
\endmatrix\right),$$
  then the intersection matrix $M_n$ of $A''$  is 
$$M_n=\pmatrix
   -\frac{n+2}{n+1}&\frac{1}{n+1}&^t\ba_1\\
   \frac{1}{n+1}&-\frac{n+2}{n+1}&^t\ba_2\\
     \ba_1& \ba_2&N
\endmatrix\right);$$ 

(ii)
  denote $K(A)=\sum m_iE_i$, and $K(A')=\sum m'_iE'_i+\sum n_jF_j$,
   then $K(A'')=\sum m'_iE''_i$;

(iii)  
  $K(A)\cdot E_i=K(A')\cdot E'_i=K(A'')\cdot E''_i$ for every 
  $i=1,\ldots , r$;
  denote this value by $c_i$;
  then $c_i\geq 0$, $c_1=c_2=0$
  and $\bm= ^t(m_1,\ldots ,m_r)$,  $\bm'= ^t(m_1',\ldots ,m_r')$
  are the solutions of linear equations $M\bm=\bc$, $M_n\bm'=\bc$,
  where  
  $\bc= ^t(c_1,\ldots ,c_r)$;

(iv)
  $-K(A')^2=-K(A'')^2=-^t\bm' M_n\bm'=-^t\bc M_n^{-1}\bc$ and 

  $-K(A)^2=-^t\bm M \bm=-^t\bc M^{-1}\bc$.
\end{lem}
\begin{pf}
  (i) follows from that $\varphi^*E''_1=E_1'+\frac{n}{n+1}F_1+\ldots
   \frac{1}{n+1}F_n$ and the similar equality for $\varphi^*E''_2$.
  As the singularity obtained by the contraction $\varphi$ is a
  rational double point, 
  it follows that $K(A')=\varphi^*K(A'')$.
  (ii) follows immediately from this fact.
  In (iii),  the equality $K(A')\cdot E'_i=K(A'')\cdot E''_i$
  follows also from this fact.
  The equality $K(A)\cdot E_i=K(A')\cdot E'_i$ follows from the fact
  that every pair $E_i$, $E'_i$ has the common arithmetic genus
  and the self intersection number by the definition of 
  (-2)-insertion.
  The statement of the value of $c_i$ follows from that 
  $A'$ has no (-1)-curve 
  and $E'_1$, $E'_2$ are 
  (-2)-curves.
  The first equality of  (iv) also follows from 
  that $K(A')=\varphi^*K(A'')$ and the rest
  is  from the definitions.
\end{pf}
\begin{lem}
\label{difference}
%**(difference)
  Let  $M_n$  be negative definite $r\times r$ real symmetric 
  matrices
  for every integer $n\geq 0$ as follows:
$$M_n=\pmatrix
   -\frac{n+2}{n+1}&\frac{1}{n+1}&^t\ba_1\\
   \frac{1}{n+1}&-\frac{n+2}{n+1}&^t\ba_2\\
     \ba_1& \ba_2&N
\endmatrix\right).$$
  Assume that every $(i,j)$-entry of $M_n$  is non-negative
  for $i\neq j$.
  Let $c_i$  $(i=1,\ldots ,r)$ be non-negative real numbers
  with $c_1=c_2=0$ and $\bc$ be $^t(c_1,\ldots c_r)$.
  Let $\bm=\ ^t(m_1,\ldots ,m_r)$ and $\bm'=\ ^t(m'_1,\ldots ,m'_r)$ 
  be the solutions of linear equations $M_0\bm=\bc$ and $M_n\bm'=\bc$,
  respectively.
  Then,

(i)
   $m_i\leq 0$ and $m'_i \leq 0$  for every $i$, and

(ii)  $$-^t\bm'M_n\bm'=-^t\bm M_0\bm+\frac{n}{n+1}(m_1-m_2)(m'_1-m'_2).$$
\end{lem}
\begin{pf}
  (i). Let $e_{ij}$  be the $(i,j)$-entry of $M_0$.
  Let $n_i$  be $m_i$  if $m_i\geq 0$ and $0$  if $m_i<0$.
  Let  $n'_i$  be $0$  if $m_i\geq 0$ and $-m_i$  if $m_i<0$.
  Then $\bn=\ ^t(n_1,\ldots , n_r)$  and  
  $\bn'=\ ^t(n'_1,\ldots , n'_r)$ satisfy that 
  $\bm=\bn-\bn'$.
  One will prove $\bn=0$.
$$0\leq\ ^t\bn\bc=\ ^t\bn M_0\bm=\ ^t\bn M_0\bn-\ ^t\bn M_0\bn'.$$
  Here $ ^t\bn M_0\bn\leq 0$ by the negative definiteness 
  and $ ^t\bn M_0\bn'\geq 0$  by $n_ie_{ij}n'_j\geq 0$ for
  every $i,\ j$.
  Therefore  $\bn$  must be $0$.
  The proof for $m'_i$  is the same.
  For the proof of (ii) let  $e_{ij}$ and $e'_{ij}$  be the $(i,j)$-entries
of $M_0$  and
  $M_n$,  respectively.
  Then 

  $^t\bm'M_n\bm'- ^t\bm M_0\bm=\sum_{i,j\geq 1}m'_im'_je'_{ij}-
  \sum_{i,j\geq 1}m_im_je_{ij}$

$  =\sum_{i\geq 1}c_i(m'_i-m_i)
  =\sum_{i\geq 3}c_i(m'_i-m_i)$

$  =\sum_{j\geq 1}m_j(\sum_{i\geq 3}m'_ie_{ij}-\sum_{i\geq 3}m_ie_{ij})$

$  = \sum_{j\geq 1}m_j(m_1e_{1j}+m_2e_{2j}-m'_1e'_{1j}-m'_2e'_{2j})$

$  =-m'_1\sum_{j\geq 1}m_je'_{1j}-m'_2\sum_{j\geq 1}m_je'_{2j}$

$  =-\sum_{i=1,2,j=1,2}m'_im_je'_{ij}
   -m'_1\sum_{j\geq 3}m_je'_{1j}-m'_2\sum_{j\geq 3}m_je'_{2j}$

$  =-\sum_{i=1,2,j=1,2}m'_im_je'_{ij}
   -m'_1\sum_{j\geq 3}m_je_{1j}-m'_2\sum_{j\geq 3}m_je_{2j}$

$  =-\sum_{i=1,2,j=1,2}(m'_im_je'_{ij}-m'_im_je_{ij})
  = -\frac{n}{n+1}(m_1-m_2)(m'_1-m'_2).$
\end{pf}

\vskip.5truecm
\begin{lem}
\label{matrix}
%**(matrix)
  Let  $\bm$, $\bm'$ and $M_n$  be as in Lemma\ref{difference}.  
  Then $$m'_1-m'_2=\frac{\det M_0}{\det M_n}(m_1-m_2).$$
\end{lem}
\begin{pf}
  In the following,
  one uses the notation in  Lemma\ref{difference}. 
  Since $\bm'$  is the solution of the linear equations:
$M_n\bm'=\bc$,   by Cramer's formula one obtains that
$$m'_1=\frac{1}{\det M_n}\left|\matrix
  0 & \frac{1}{n+1} & ^t\ba_1\\
  0 & -\frac{n+2}{n+1} & ^t\ba_2\\
  \bc' & \ba_2 & N
\endmatrix\right|$$
$$m'_2=\frac{1}{\det M_n}\left|\matrix
  -\frac{n+2}{n+1} & 0 & ^t\ba_1\\
  \frac{1}{n+1} & 0 &  ^t\ba_2\\
  \ba_1 & \bc' & N
\endmatrix\right|$$
where $\bc'=^t(c_3,\ldots ,c_r)$.
  For $i=1,2$, let $a_i$  be
$$\left|\matrix
  0 & ^t\ba_i \\
  \bc' & N
\endmatrix\right|,$$
  then
$$m_1'=\frac{1}{\det M_n}\{-\frac{1}{n+1}a_2-\frac{n+2}{n+1}a_1
+({\rm terms\ independent\ of}\ n)\}$$
$$m'_2=\frac{1}{\det M_n}\{-\frac{n+2}{n+1}a_2-\frac{1}{n+1}a_1
+( {\rm terms\ independent\ of}\ n)\}.$$
  Therefore
$$m'_1-m'_2=\frac{1}{\det M_n}\{a_2-a_1+({\rm terms\ independent\ of}\ n)\}.$$
  Now comparing this with the equality for $n=0$,
  it follows that 
$$\det M_n(m'_1-m'_2)=\det M_0(m_1-m_2).$$
\end{pf}
\vskip.5truecm
\begin{cor}
\label{positive}
  Let $A'=\sum_{i=1}^r E'_i+\sum_{j=1}^nF_j$  be 
  a (-2)-insertion of 
  $A=\sum _{i=1}^rE_i$.
  Then 
$$-K(A')^2\geq -K(A)^2.$$
\end{cor}
\begin{pf}
  Under the notation  in \ref{property},
  $-K(A')^2=- ^t{\bm}'M_n{\bm}'$ and $ -K(A)^2=- ^t{\bm}M_0{\bm}$.
  Then by \ref{difference} and \ref{matrix} it follows the 
  equality 
  $$-K(A')^2=-K(A)^2+\frac{n\det M_0}{(n+1)\det M_n}(m_1-m_2)^2.$$
  Here noting that $M_0$ and $M_n$ are negative definite matrices
  of the same size,
  one obtains the inequality ${\det M_0}/{\det M_n}>0$,
  which completes the proof.
\end{pf}
\begin{thm}
\label{DCC}
%**(DCC)
  The set $\{-K^2| K $ is the numerical canonical divisor of
   a normal surface singularity $(X,x)\}$
  has no accumulation points from above.
\end{thm}
\begin{pf}
  Assume there exists an accumulation point $\alpha$ from above.
  Then one can take a sequence of surface singularities 
  $\{(X_{(n)},x_{(n)})\}_{n\in \bN}$,
  such that $\alpha = \lim_{n\to \infty}(-K_{(n)}^2)$ and
  $-K_{(n)}^2>-K_{(n+1)}^2$ for  all $ n$,
  where  $K_{(n)}$ is the numerical canonical divisor of 
  the singularity 
  $(X_{(n)},x_{(n)})$.
  Let $A_{(n)}$  be the exceptional divisor on the minimal resolution
  of $(X_{(n)},x_{(n)})$, $E_{(n)i}$ and $F_{(n)j}$ be non-(-2)-curve
  and (-2)-curve of  $A_{(n)}$, respectively.
  Since  $-K_{(n)}^2$  are bounded and   
  $(K_{(n)}\cdot E_{(n)i})^2/(-E_{(n)i}^2)\geq 1/3$, 
  $K_{(n)}\cdot F_{(n)j}=0$, 
  the numbers of the components $E_{(n)i}$
  are bounded by \ref{basic lemma}.
  Replacing by a subsequence,
  one may assume that the numbers of components $E_{(n)i}$'s 
  are the same;
  say $r$. 
  By \ref{basic lemma} and the equality   
$$(K_{(n)}\cdot E_{(n)i})^2/(-E_{(n)i}^2)=
(2p_a(E_{(n)i})-2)^2/(-E_{(n)i}^2)+2(2p_a(E_{(n)i})-2)-E_{(n)i}^2,$$ 
  it follows that $p(E_{(n)i})$  and $-E_{(n)i}^2$ are bounded;
  therefore one may assume that these are constant by taking
  a subsequence.
  Since $(E_{(n)i}\cdot E_{(n)j})^2<E_{(n)i}^2 E_{(n)j}^2$,
  the intersections $E_{(n)i}\cdot E_{(n)j}$  are also bounded,
  so one may assume that these are also constant for every $(i,j)$.
 Now one may assume that the numerical conditions 
 of the non-(-2)-curves $\{E_{(n)i}\}_{i=1}^r$ are the same for every $n$.
  Next one should note that the numbers of the connected components
  of $F_{(n)j}$'s are bounded.
  In fact, 
  every connected component has an intersection with at least one of
  the $E_{(n)i}$'s because of the connectedness of the exceptional
  divisor,
  and the numbers of the connected components
  of $F_{(n)j}$'s which  $E_{(n)i}$ can intersect
  are bounded
  because of the negative definiteness of the intersection matrix.
  So one may assume that the numbers of the connected components of 
  $F_{(n)j}$'s are constant; say $s$.
   Let  $S_{(n)k}$, $k=1,\ldots , s$  be the connected components
  of $F_{(n)j}$'s.
   As the intersections of $S_{(n)k}$'s and $E_{(n)i}$'s are
   bounded, 
  one may assume that each $S_{(n)k}$ intersects the same 
  $E_{(n)i}$'s with the same intersection numbers for all  $n$.
  If the numbers of irreducible components of $S_{(n)k}$
  are bounded for all $n$,
  then one may assume that  the configuration of 
  $S_{(n)k}\cup (\bigcup_{i=1}^rE_{(n)i})$ 
  is fixed.
  Note that $S_{(n)k}$ is one of the configurations A$_m$
  $(m\geq 1)$, D$_m$ $(m\geq 4)$, E$_6$, E$_7$, E$_8$.
  If  the number of  irreducible components
 of $S_{(n)k}$  tends to $\infty$,
  then $S_{(n)k}$ is either A$_m$ or  D$_m$.
  Therefore
  one may assume that $A_{(n)} \to A_{(n+1)}$  is a composite
  of (-2)-insertions.
  By \ref{positive} 
  one obtains that
$$-K_{(n)}^2\leq -K_{(n+1)}^2,$$
  which is a contradiction.
\end{pf}
\vskip.5truecm
\section{\bf Accumulation points of $-K^2$}
\begin{say}
  Let  $A$  be a compact divisor on a non-singular surface.
  A string $S$  in $A$ is a chain of (-2)-curves
  $A_1,\ldots , A_n$ so that $A_i\cdot A_{i+1}=1$ 
  $(i=1, \ldots , n-1)$,
  and these account for all intersections in $A$  except that
  at least one of $A_1$ and $A_n$ intersects other curves.
\end{say}
\vskip.5truecm
  Now one prepare an easy lemma for (-2)-insertions. 

\begin{lem}
\label{maxstring}
  Let $A=\sum _{i=1}^rE_i \to  A'=\sum_{i=1}^r E'_i+S_1+S_2$ 
  be the composite of (-2)-insertions with the insertion strings
  $S_1$ and $S_2$.
  If $S_1$  and $S_2$ are in one maximal string of $A'$,
  then $A\to A'$  is regarded as one (-2)-insertion. 
\end{lem}

\begin{thm}
\label{rational}
%**(rational)
  All accumulation points of the set 
  $\{-K^2| K $ is the numerical canonical divisor of
  a normal surface singularity $(X,x)\}$
  from below are rational numbers.
\end{thm}
\begin{pf}
  Let  $\alpha$  be an accumulation point of the set $\{-K^2\}$
  from below.
  Then one can take a sequence of surface singularities 
  $\{(X_{(n)},x_{(n)})\}_{n\in \bN}$ such that 
  $\alpha = \lim_{n\to \infty}(-K_{(n)}^2)$ and
  $-K_{(n)}^2<-K_{(n+1)}^2$ for every $ n$,
  where  $K_{(n)}$ is the numerical canonical divisor of 
  $(X_{(n)},x_{(n)})$.
  Let $A_{(n)}$  be the exceptional divisor on the minimal resolution
  of $(X_{(n)},x_{(n)})$
  Replacing by a suitable subsequence,
  one may  assume that $A_{(1)}\to A_{(n)}$ is a composite of
  (-2)-insertions:
  $A_{(1)}=A_{0\ldots 0}=\sum _{i=1}^rE_i$,
  $A_{(n)}=A_{n_1\ldots n_k}=
  \sum_{i=1}^rE'_i+\sum_{j=1}^kS_j$,
  where $S_j=\sum_{s=1}^{n_j}F_{js}$  is the insertion string
  with the props $E'_{2j-1}$, $E'_{2j}$. 
  And one may assume that if $n\to \infty$,
  then $n_j\to \infty$  for all $j=1,\ldots , k$.
  By \ref{maxstring}, we may assume that every maximal string
  of $A_{(n)}$ has at most one insertion string. 
  Let  $\varphi$  be the contraction morphism of $S_1,\ldots , S_k$.
  Denote $\varphi_*A_{n_1\ldots n_k}$ by $A''_{n_1\ldots n_k}$ and
  $\varphi_*E'_i$  by $E''_i$.
  Let $M_{n_1\ldots n_k}$ be the intersection matrix of $A''_{n_1\ldots n_k}$,
  then by (i) of \ref{property}
$$M_{n_1\ldots n_k}=\pmatrix
   -\frac{n_1+2}{n_1+1}&\frac{1}{n_1+1}&\ldots &\ldots &0 &0&^t\ba_1\\
   \frac{1}{n_1+1}&-\frac{n_1+2}{n_1+1}& \ldots&\ldots& 0&0&^t\ba_2\\
  \vdots&\vdots&\ddots&  &\vdots &\vdots &\vdots\\
  \vdots&\vdots&  &\ddots&\vdots &\vdots &\vdots\\
  0&0&\ldots &\ldots & -\frac{n_k+2}{n_k+1}&\frac{1}{n_k+1}&^t\ba_{2k-1}\\
  0&0&\ldots &\ldots & \frac{1}{n_k+1}&-\frac{n_k+2}{n_k+1}&^t\ba_{2k}\\
     \ba_1& \ba_2 &\ldots &\ldots &\ba_{2k-1} &\ba_{2k} &N
\endpmatrix.$$
  Denote $\lim_{n_1,\ldots ,n_k\to \infty}M_{n_1\ldots n_k}$
  by $M_{\infty\ldots \infty}$.
  Then by the definition of $M_{n_1\ldots n_k}$,
  every entry of 
  the matrix $M_{\infty\ldots \infty}$ is rational.
  By \ref{property}, $K(A''_{n_1\ldots n_k})\cdot E''_i
  =K(A_{0\ldots 0})\cdot E_i$  for every $i$ and $n_1,\ldots , n_k$.
  Denote this value by $c_i$, and put $\bc=\  ^t(c_1,\ldots ,c_r)$,
  then $c_i$ is non-negative integer for every 
  $i$ and $c_i=0$ for $i=1,\ldots ,2k$.
  Noting that $-K(A_{n_1 \ldots n_k})^2=-{K(A''_{n_1 \ldots n_k})}^2=
  - ^t\bc M_{n_1 \ldots n_k}^{-1} \bc$,
  if one proves that $M_{\infty\ldots \infty}$ is invertible,
  it follows that 
  $\lim_{n_1,\ldots , n_k\to \infty}-K(A_{n_1 \ldots n_k})^2=
  -  ^t\bc M_{\infty \ldots \infty}^{-1} \bc$,
  which is rational as required.
  First one claims that $M_{\infty n_2 \ldots n_k}$  is invertible
  for every $(n_2,\ldots ,n_k)$.
  For this,
  it is sufficient to prove that,
  fixing $(n_2,\ldots , n_k)$ there exists a positive number 
  $\epsilon$ such that 
  $|\det M_{n_1 \ldots n_k}|>\epsilon$ for all $n_1$.
  By the proof of 1.9,
$$-K(A_{n_1n_2\ldots n_k})^2=-K(A_{0n_2\ldots n_k})^2
  +\frac{n_1 \det M_{0n_2\ldots n_k}}
  {(n_1+1)\det M_{n_1n_2\ldots n_k}}(m_1-m_2)^2.$$
  This shows that $m_1\neq m_2$,
  because $-K^2$  is strictly increasing.
  Since  $-K(A_{n_1n_2\ldots n_k})^2<\alpha$,
  it follows that 
$$|\det M_{n_1n_2\ldots n_k}|=
\frac{n_1|\det M_{0n_2\ldots n_k}|(m_1-m_2)^2}
      {(n_1+1)(-K(A_{n_1n_2\ldots n_k})^2+K(A_{0n_2\ldots n_k})^2)}$$
$$> \frac{|\det M_{0n_2\ldots n_k}|(m_1-m_2)^2}
        {2(\alpha +K(A_{0n_2\ldots n_k})^2)}
>0$$
  for all $n_1$.

  Next fixing $(n_3,\ldots, n_k)$,
  one will show that $M_{\infty \infty n_3 \ldots n_k}$ 
  is invertible.
  Take the solution  $\bm=\ ^t(m_1,\ldots ,m_r)$ 
  for linear equations 
  $M_{\infty 0 n_3\ldots n_k}\bm=\bc$.
  Then again by the proof of \ref{positive},
$$ ^t\bc M^{-1}_{\infty n_2\ldots n_k}\bc
=  ^t\bc M^{-1}_{\infty 0n_3\ldots n_k}\bc 
  + \frac{n_2\det M_{\infty 0n_3\ldots n_k}}
         {(n_2+1)\det M_{\infty n_2\ldots n_k}}
   (m_3-m_4)^2.$$
  If $m_3\neq m_4$,
  then one can proceed to get the regularity of 
  $M_{\infty\infty n_3\ldots n_k}$ 
  in the same way as the discussion on $n_1$.
  If $m_3=m_4$,
  then $^t\bc M^{-1}_{\infty n_2\ldots n_k}\bc$  is 
  constant for every  $n_2$.
  In this case we define $M_{\infty\infty n_3\ldots n_k}$
  to be  $M_{\infty 0 n_3\ldots n_k}$
  which is invertible.    

  By the successive procedure,
   it finally follows that  $M_{\infty\ldots \infty}$  is invertible.
\end{pf} 
\vskip.5truecm
\begin{prop}
    For an arbitrary rational number  $r$,
  there exists an integer $n$  such that 
   $n+r$  is an accumulation point of $\{-K^2\}$.
\end{prop}
\begin{pf}
  One may assume that  $r=k/m$ with $k,m\in \bN$  and
  $k$ is even.
  Let $A_{(n)}$ be a divisor $E_1+E_2+\sum_{j=1}^nF_j$
  with the following dual graph:

$\ \ \ \ \ x-o-\cdots - o-o-o\ ,$

$\ \ \ \ \ \ \ \ \ \ \  |$

$\ \ \ \ \ \ \ \ \ \ \  x$
\newline
where $x=E_i$ $(i=1,2)$ and $o=F_j$ $(j=1,\ldots , n)$.
  Assume that $F_j $'s are (-2)-curves, and $E_i$'s are non-singular
  curves of genus $k/2$ with self-intersection number $-2km-2$.
  First note that these divisors can be constructed in ruled
  surfaces of genus $k/2$.
  Then one can show that the accumulation point of the set
  $\{-K(A_{(n)})^2\}$ is  $k/m$+integer.
  In fact,
  put $K(A_{(n)})=\sum_{i=1}^2m_{(n)i}E_i+\sum_{j=1}^ns_{(n)j}F_j$,
  then, solving the linear equations,
  one obtains:
  $$m_{(n)1}=m_{(n)2}=\frac{(1+1/n)k(2m+1)}{(2km+2)(-1-1/n)+2}.$$
Therefore
$$-K(A_{(n)})^2=-\sum_{i=1}^2m_{(n)i}K(A_{(n)})\cdot E_i=
\frac{2(1+1/n)k^2(2m+1)^2}{(2km+2)(1+1/n)-2}$$
which tends to 
$\frac{k^2(2m+1)^2}{km}=k/m$+integer,
when $n$ tends to $\infty$.
\end{pf}
\vskip.5truecm
\begin{prop}
\label{r,k}
  For positive integers $r$, $k$ with $r\geq k-1$,
  the rational number  $r^2/k$  is an accumulation point
  of $\{-K^2\}$.
\end{prop}
\begin{pf}
  For $r$, $k$  such that $r-k\equiv 1 (\text{mod} 2)$,
  take a divisor with the following configuration:

$o-\cdots -o-x$ ,
\newline
where $o$'s are (-2)-curves and $x$  is a non-singular curve
  of genus  $(r-k+1)/2$ with the self-intersection number 
  $-(k+1)$.
  By the condition of $r$, $k$, 
  the curve $x$  is not a (-1)-curve;
  therefore this divisor can be the exceptional divisor 
  of the minimal resolution of a singularity.
  Let the number of (-2)-curves be $n$,
  then, for the numerical canonical divisor $K_{(n)}$ of
  the corresponding singularity,
  it follows that 
$$-K_{(n)}^2=\frac{r^2}{k+1-n/(n+1)}.$$
Therefore if $n$ tends to $\infty$, then this value tends to
  $r^2/k$.

  For $r$, $k$ with $r-k\equiv 0 (\text{mod} 2)$,
  take a divisor with the following configuration:

$o-\cdots -o-x-o-\cdots -o$,

where $o$'s are (-2)-curves and $x$ is non-singular curve
  of genus $(r-k)/2$ with the self-intersection number
  $-(k+2)$.
  Then this divisor also can be the exceptional divisor 
  of the minimal resolution of a singularity.
  Let the numbers of (-2)-curves of the both hand sides of $x$
  be 
  $n$  and $s$.
  Then 
$$-K_{(n,s)}^2=\frac{r^2}{k+2-n/(n+1)-s/(s+1)}.$$
  Therefore if $n$ and $s$ tend to $\infty$, 
  then this value tends to
  $r^2/k$.
\end{pf}
\begin{cor}
  An arbirary positive integer is an accumulation point
  of $\{-K^2\}$.
\end{cor}
\begin{pf}
In Proposition\ref{r,k},
  put $r=k$.
\end{pf}   
\vskip.5truecm
\begin{exmp}
  Here one calculates the value of $-K^2$ for rational triple points.
  Rational triple points are classified into 9 classes
  according to the dual graphs in \cite{Artin}.
  In the following dual graphs,
  $x$ denotes the (-3)-curve and $o$ denotes (-2)-curve.

$I_{n,s,t}: o-\cdots -o-x-o-\cdots -o$

$\ \ \ \ \ \ \ \ \ \ \ \ \ \ \ \ \ \ \ \ \ \ \ \ \ \ |$

$\ \ \ \ \ \ \ \ \ \ \ \ \ \ \ \ \ \ \ \ \ \ \ \ \ \ o-\cdots -o$  

where $n$, $s$, $t$ are the numbers of (-2)-curves in each string.
  Then 
$$-K(I_{n,s,t})^2=\frac{1}{3-n/(n+1)-s/(s+1)-t/(t+1)}.$$

$II_{n,s}:  o-\cdots-o-x-o-\cdots -o-o$

$\ \ \ \ \ \ \ \ \ \ \ \ \ \ \ \ \ \ \ \ \ \ \ \ \ \ \ \ \ \ \ \ \ \ \ \ \
\ \ \ \ \ \ \ \ |$

$\ \ \ \ \ \ \ \ \ \ \ \ \ \ \ \ \ \ \ \ \ \ \ \ \ \ \ \ \ \ \ \ \ \ \ \ \
\ \ \ \ \ \ \ \ o$

where $n$ is the number of (-2)-curves in the left string and $s$
  is that in the graph of the right hand side of $x$.
  Then 
$$-K(II_{n,s})^2=\frac{n+1}{n+2}.$$

$III_{n,s}: o-\cdots -o-x-o-\cdots -o$

$\ \ \ \ \ \ \ \ \ \ \ \ \ \ \ \ \ \ \ \ \ \ \ \ \ \ \ \ \ \ \ \ \ \ |$

$\ \ \ \ \ \ \ \ \ \ \ \ \ \ \ \ \ \ \ \ \ \ \ \ \ \ \ \ \ \ \ \ \ \ o$

where  $n$ is the number of (-2)-curves in the left string and $s$
  is that in the graph of the right hand side of $x$.
  Then 
$$-K(III_{n,s})^2=\frac{1}{3-n/(n+1)-2(s-1)/(s+1)}.$$

$IV_{n}: o-\cdots -o-x-o-o-o-o$

$\ \ \ \ \ \ \ \ \ \ \ \ \ \ \ \ \ \ \ \ \ \ \ \ \ \ \ \ \ \ \ \ \ \ \ \ |$

$\ \ \ \ \ \ \ \ \ \ \ \ \ \ \ \ \ \ \ \ \ \ \ \ \ \ \ \ \ \ \ \ \ \ \ \ o$

where $n$ is the number of (-2)-curves in the left string.
Then
$$-K(IV_n)^2=\frac{4n+4}{3n+7}.$$

$V_n: o-\cdots-o-x-o-o-o-o-o$

$\ \ \ \ \ \ \ \ \ \ \ \ \ \ \ \ \ \ \ \ \ \ \ \ \ \ \ \ \ \ \ \ \ \ \ \ \
\ \ \ |$

$\ \ \ \ \ \ \ \ \ \ \ \ \ \ \ \ \ \ \ \ \ \ \ \ \ \ \ \ \ \ \ \ \ \ \ \ \
\ \ \ o$

where $n$ is the number of (-2)-curves in the left string.
Then
$$-K(V_n)^2=\frac{3n+3}{2n+5}.$$

$VI_n: o-\cdots -o-o-o-o$ 

$\ \ \ \ \ \ \ \ \ \ \ \ \ \ \ \ \ \ \ \ \ \ \ \ \ \ |$

$\ \ \ \ \ \ \ \ \ \ \ \ \ \ \ \ \ \ \ \ \ \ \ \ \ \ x$

where $n+2$ is the number of (-2)-curves.
Then
$-K(VI_n)^2=(n+3)/9.$
\vskip.5truecm 

$VII:x-o-o-o-o-o$  Then $-K(VII)^2=2/3.$

$\ \ \ \ \ \ \ \ \ \ \ \ \ \ \ \ \ \ \ |$

$\ \ \ \ \ \ \ \ \ \ \ \ \ \ \ \ \ \ \ o$
\vskip.5truecm

$VIII:x-o-o-o-o-o-o$   Then  $-K(VIII)^2=4/5$.

$\ \ \ \ \ \ \ \ \ \ \ \ \ \ \ \ \ \ \ \ \  |$

$\ \ \ \ \ \ \ \ \ \ \ \ \ \ \ \ \ \ \ \ \  o$
\vskip2truecm

$IX:x-o-o-o-o-o-o$ Then $-K(IX)^2=2/3$.

$\ \ \ \ \ \ \ \ \ \ \ \ \ \ \ \ \ \ \ \ \ \ \ \ \ \ \ \ |$

$\ \ \ \ \ \ \ \ \ \ \ \ \ \ \ \ \ \ \ \ \ \ \ \ \ \ \ \ o$
\end{exmp}
\vskip.5truecm
\begin{prop}
\label{interval}
  The accumulation points of $-K^2$ for normal surface singularities
  in the interval $[0, 1]$ are
  $m/(m+1)$ $(m\in \bN) $  and 1.
\end{prop}
\begin{pf}
  Let  $\alpha\in [0,1]$  be an accumulation point of $-K^2$.
  Then one can take a sequence of normal surface singularities
  $\{(X_{(n)},x_{(n)})\}_{n\in \bN}$ whose numerical canonical
  divisors  $K_{(n)}$ satisfy $0<-K_{(n)}^2< 1$ and 
  $\lim_{n\to \infty}(-K_{(n)}^2)=\alpha$.
  Then $(X_{(n)},x_{(n)})$ are rational triple  points.
  In fact,
  for the fundamental cycle $Z$ on the minimal resolution of 
  $(X_{(n)},x_{(n)})$,
  the quadratic form $(xK+yZ)^2$  is negative semidefinite,
  where $K=K_{(n)}$.
  Therefore  $K^2Z^2\geq (K\cdot Z)^2$.
  Since $1>-K^2$ and $K\cdot Z$  is a non-negative integer,
  one obtains the inequalities:
$$(\ref{interval}.1)\ \ \ \ \ \ \ \ \ \ \ \ \ \ \ \ 
  -Z^2> (-K^2)(-Z^2)\geq (K\cdot Z)^2 \geq K\cdot Z.$$
  Hence  $Z^2+K\cdot Z<0$,
  which yields that $Z^2+K\cdot Z=-2$
  because $p_a(Z)=1+\frac{1}{2}(Z^2+K\cdot~Z)\geq 0$.
  Thus $(X_{(n)}, x_{(n)})$  is a rational singularity.
  On the other hand, $-Z^2>(K\cdot Z)^2$ by (\ref{interval}.1).
  This gives the inequality:
$$0>Z^2+(K\cdot Z)^2=(K\cdot Z)^2-K\cdot Z -2= 
  (K\cdot Z+1) (K\cdot Z-2)$$
  by combining $Z^2=-K\cdot Z -2$.
  Hence  $K\cdot Z= 0$ or 1,
  equivalently $Z^2=-2 $ or $-3$,
  which shows that the multiplicity of $(X_{(n)},x_{(n)})$
  is 3 since $-K^2>0$.

  For  rational triple points,
  we have all accumulation points in the previous example.
\end{pf}
%%%%%%%%%%%%%%%%%%%%%%%%  3  %%%%%%%%%%%%%%%%%%%%%%%%%
\vskip.5truecm
\section{\bf Boundedness of invariants for the constant $-K^2$}

  In this section one  observes the behavior of various invariants
  of a singularity under fixing $-K^2$.
\begin{exmp}
  Under fixing $-K^2$, the numbers of the exceptional curves on
  minimal resolutions are not in general bounded.
  For example A$_n$ has $n$-exceptional curves and $-K^2=0$ 
  for every $n$.

  One can also see such an example with $-K^2\neq 0$.
  Take singularities with the minimal 
  resolutions of the following graphs:
  $$x-o-\cdots -o-x,$$ 
   where $x$'s are  (-3)-curves and 
  $o$'s are  
   (-2)-curves and the number of (-2)-curves is $n$.
  Then $-K^2=1$  for every $n$.
\end{exmp}
  
  By these examples one  can see that the boundedness,
  as  in \cite{Al}, for singularities under fixing $-K^2$  cannot be
  expected.
  But still $-K^2$ has a power to controll other invariants of 
  singularities.

\begin{thm}
\label{mult}
%**(mult)
  For a normal surface singularity $(X,x)$,
  the numerical canonical divisor  $K$, the multiplicity 
   and the embedding dimension of $(X,x)$ 
  satisfy the following inequalities:
\end{thm}
$$ -K^2\geq{\text{{mult}}} (X,x) -4,$$
$$  -K^2\geq{\text{{embdim}}}(X,x) -5.$$

\begin{pf}
  Let  $f:Y\to X$  be the minimal resolution of the singularity
  $(X,x)$ and   
  $K$  the numerical canonical divisor.   
  Let $h:Y'\to Y$ be a suitable proper birational morphism
  such that $g:=f\circ h:Y'\to X$  factors through the blowing-up
  of the maximal ideal  ${\frak m}_{X,x}$.
  Let  $D$  be the divisor such that 
  ${\cal O}_{Y'}(-D)={\frak m}_{X,x}{\cal O}_{Y'}$.
  Then $-D$  is relatively nef with respect to $g$ and by 
  \cite{To},
  3.12, it follows that  $p_a(D)\geq 0$.
  If one denote by $K'$ the numerical canonical divisor for 
  $g:Y'\to X$, it follows that $h^*K=K'-\Delta$ with $\Delta>0$.
  Noting that $(h^*K+D)^2\leq 0$, 
  one obtains

$$-K^2=-(h^*K)^2\geq 2(h^*K\cdot D+D^2)-D^2$$

$$\ \ \ \geq 2(K'\cdot D+D^2)-D^2\geq -4-D^2
  =-4+{\text{mult}} (X,x).$$

  On the other hand,
  by \cite{Ab} ${\text{embdim}}R\leq {\text{mult}}_{\frak m}R+\dim R -1$
  for a Cohen-Macaulay local ring $(R,{\frak m})$.
  Therefore  follows the inequality involving the embedding dimension.
\end{pf}

Immediately by this theorem one obtains
the following corollary.

\begin{cor}
   The multiplicities and the embedding dimensions are bounded
  for the constant $-K^2$.
\end{cor}

\begin{thm}
\label{agenus}
%***(agenus)
  For a normal surface singularity $(X,x)$,
  the numerical canonical divisor and the arithmetic genus 
  satisfy:
$$-K^2\geq 4p_a(X,x)-3.$$
\end{thm}

\begin{cor}
  The arithmetic genera are bounded for the constant $-K^2$.
\end{cor}
\begin{pf}
  By the negative definiteness,
  $(K+D)^2\leq 0$  for every divisor $D$ with the support on the exceptional
  divisor.
  Therefore $-K^2\geq -D^2+2(K\cdot D+D^2)=-D^2+4(p_a(D)-1)$,
  where $-D^2\geq 1$.
\end{pf}

\begin{thm}
\label{nindex}
%***(nindex)
   The numerical indices are bounded for the constant $-K^2$.
\end{thm}
\begin{pf}
  Assume that the numerical index $r(X,x)$  is not bounded
  for the constant $-K^2$.
  Then there exists a sequence of singularities 
  $\{(X_{(n)},x_{(n)})\}_{n\in \bN}$ with constant $-K^2$  
  such that  $r(X_{(n)},x_{(n)})<r(X_{(n+1)},x_{(n+1)})$
  for every $n$.
  One will show a contradiction 
  by constructing a subsequence with the constant $r(X_{(n)},x_{(n)})$.
  Let  $A_{(n)}$  be the exceptional divisor on the minimal resolution
  of   $(X_{(n)},x_{(n)})$.
  Then,
  replacing by a suitable subsequence,
  one may assume that  $A_{(n)}\to A_{(n+1)}$  is a composite of 
  (-2)-insertions.
  Therefore it is sufficient to prove 
  that a (-2)-insertion 
  $A=\sum_{i=1}^rE_i \ \to A'=\sum_{i=1}^rE'_i+\sum_{j=1}^nF_j$
  with $-K(A)^2=-K(A')^2$ does not change the numerical indices.
  Denote $K(A)$ and $K(A')$  by $\sum_{i=1}^rm_iE_i$  and
  $\sum_{i=1}^rm'_iE'_i+\sum_{j=1}^nn_jF_j$.
  By \ref{property}, \ref{difference} and \ref{matrix},
  $$-K(A')^2= -K(A)^2+\frac{n}{n+1}\frac{\det M}{\det M_n}(m_1-m_2)^2.$$
  Then $m_1=m_2$ because $-K(A')^2= -K(A)^2$.
  Let $K'$  be a $\bQ$-divisor 
  $\sum_{i=1}^rm_iE'_i+m_1\sum_{j=1}^nF_j$,
  then it satisfies the equations 
  $K'\cdot E'_i=2p_a(E_i)-2-E_i^2=K(A')\cdot E'_i$
  and $K'\cdot F_j=2p_a(F_j)-2-F_j^2=K(A')\cdot F_j$  for all $i$ and $j$.
  By the uniqueness of the solution of these equations,
  it follows that $K'=K(A')$.
  Hence the sets $\{m'_i, n_j\}$  and $\{m_i\}$  coincide,
  showing that the numerical indices are the same.
\end{pf}

  Immediately by this theorem one obtains the following corollary.

\begin{cor}
  For rational singularities,
  the indices are bounded for the fixed $-K^2$.
\end{cor}

\begin{thm}
\label{ggenus}
%**(ggenus)
  If $-K^2$  is fixed, then the geometric genera  $p_g$
  are bounded and also the plurigenera  $\gamma_m$  are bounded
  for every $m\in \bN$.
\end{thm}
\begin{pf}
  Let $r$  be the numerical
  index of the singularity $(X,x)$.
  By Kato's inequality ((16) of \cite{Ka}) 
    $p_g(X,x)\leq r(-K^2)$ and 
  by \ref{nindex}, the boundedness of the geometic genera 
  follows.
  For the plurigenera,
  recall the formula in \cite{Mo}:
$$\gamma_m(X,x)=\frac{-K^2}{2}m(m-1)+ p_g(X,x) + \epsilon,$$
  where $\epsilon$  is bounded.
  Hence the boundedness of the plurigenera also follows.
\end{pf}

  Recently Tomari obtained an inequality that
   $p_g$ is bounded by $(constant)(-K^2)$
  from above without the numerical index(\cite{To2}).

\begin{say}
  The converse of \ref{mult}, \ref{agenus},  \ref{nindex}
 and of \ref{ggenus}  do not hold.
  In fact, 
  for rational triple point $(X,x)$,
  obviously ${\text{mult}} (X,x)=3$ and ${\text{embdim}}(X,x)=4$,
  but $-K(I_{n,s,t})^2\to \infty$ as $n,s,t \to \infty$.
  As another example, let  $(X_{(n)},x_{(n)})$  be a simple elliptic singularity
  such that the exceptional divisor $E_{(n)}$ on the minimal resolution\  
  has the self intersection number  $E_{(n)}^2=-n$.
  Then $r(X_{(n)},x_{(n)})= 1 $ and 
  $p_g(X_{(n)},x_{(n)})=p_a(X_{(n)},x_{(n)})=1$ for
  every $n$,
  but $-K_{(n)}^2\to \infty$  as $n\to \infty$.
\end{say}

\begin{say}
\label{quotient}
%**(quotient)
  One may expect that $-(K+\Delta)^2$ also has D.C.C. for 
  relatively nef  $K+\Delta$, 
  where $\Delta=\sum b_iA_i$  
  with all coefficients $b_i$ belonging to a D.C.C. set ${\cal C}$,
  as is conjectured in \cite{Ko} for
  the case of compact surfaces 
  and answered in \cite{Al} for that case.
  For a singularity case,
  as  in this paper, it is true for ${\cal C}=\{0\}$,
  but not true in general.
  Indeed, take  ${\cal C}=\{1\}$.
  Let $(X,x)$  be a non-log-canonical  singularity with 
$\bC^*$-action.
  For example, take a singularity with $\bC^*$-action with the 
  following dual graph:
\newpage
$\ \ \ \ \ \ y$

$\ \ \ \ \ \ |$

$y- z-y$ ,

$\ \ \ \ \ \ |$

$\ \ \ \ \ \ y$
\newline
where  $y$'s are (-3)-curves and $z$  is a (-5)-curve.
  Then $(X,x)$  is a non-log-canonical singularity.
  Let  $(X_{(n)},x_{(n)})$  be the quotient of $(X,x)$
  by the subgroup
  $<\epsilon_n>\subset \bC^*$,
  where $\epsilon_n$  is a primitive $n$-th root of unity.
  Let  $\varphi^*_{(n)}:Y_{(n)}\to X_{(n)}$  be the log-canonical model;
  i.e. $(Y_{(n)}, A_{(n)})$  has at worst log-canonical
   singularities and 
  $K_{(n)}+A_{(n)}$  is  $\varphi_{(n)}$-ample,
  where $K_{(n)}$  is the numerical canonical divisor on $Y_{(n)}$
  and $A_{(n)}$  is the reduced exceptional divisor.
  Here note that 
  $-(K_{(n)}+A_{(n)})^2$ is $-P^2$  of the singularity 
  $(X_{(n)},x_{(n)})$.
  Then by \cite{Wah} $-(K_{(1)}+A_{(1)})^2\geq -n(K_{(n)}+A_{(n)})^2$.
  Therefore $-(K_{(n)}+A_{(n)})^2 \to 0$ as $n\to \infty$.
\end{say}

%%%%%%%%%%%%%%%%%%%Reference%%%%%%%%%%%%%%%%%%%%

%\newpage
\makeatletter \renewcommand{\@biblabel}[1]{\hfill#1.}\makeatother

\end{document}